\documentclass[a4paper,12pt]{article}

\usepackage[
  a4paper,
  top=2.0cm,
  bottom=2.0cm,
  left=2.0cm,
  right=2.0cm,
  headheight=15pt,
  headsep=0.8cm,
  footskip=1.2cm
]{geometry}
\usepackage{setspace}
\usepackage{amsmath,amssymb,amsfonts}
\usepackage{amsthm}
\usepackage{thmtools}
\usepackage{mathtools}
\usepackage{bm}
\usepackage[title]{appendix}
\usepackage{textcomp}
\usepackage{manyfoot}
\usepackage{booktabs}
\usepackage{listings}
\usepackage{hyperref}
\usepackage[capitalise,nameinlink,noabbrev]{cleveref}
\usepackage{xcolor}
\usepackage{graphicx}
\usepackage{float}
\usepackage{indentfirst}
\usepackage{enumerate}
\usepackage{stmaryrd}
\usepackage{multicol}
\usepackage{multirow}
\usepackage{array}
\usepackage{subcaption}
\usepackage{rotating}
\usepackage{microtype}
\usepackage{algorithm}
\usepackage{algorithmicx}
\usepackage{algpseudocode}
\usepackage{longtable}
\usepackage[normalem]{ulem}
\usepackage{pifont}
\usepackage[numbers]{natbib}
\graphicspath{{figures/}}

\hypersetup{
  colorlinks=false,
  hidelinks,
  hypertexnames=false
}

\theoremstyle{plain}

\theoremstyle{definition}

\theoremstyle{remark}

\newcommand{\TITLE}{A hybrid combinatorial-continuous strategy for solving molecular distance geometry problems}

\newcommand{\FUNDING}{This research was partially funded by the Brazilian research agencies FAPESP and CNPq, as well as the French research agency ANR.}
\title{\TITLE\thanks{\FUNDING}}
\author{Leonardo D. Secchin \and Wagner da Rocha \and Mariana da Rosa \and Leo Liberti \and Carlile Lavor}
\date{\today}

\begin{document}
\maketitle

	\begin{abstract}
		The Molecular Distance Geometry Problem (MDGP) is essential to structural biology, as it seeks to determine three-dimensional protein structures from partial interatomic distances. Its discretizable subclass (DMDGP) admits an exact combinatorial formulation that enables efficient exploration of the search space. However, in practical settings such as Nuclear Magnetic Resonance (NMR) spectroscopy, distances are available only within uncertainty bounds, leading to the interval variant (\emph{i}DMDGP). We propose a hybrid combinatorial--continuous framework for solving the \emph{i}DMDGP. The method combines an enumeration process derived from the DMDGP with a continuous refinement stage that minimizes a nonconvex stress function penalizing  deviations from admissible distance intervals. This integration supports a systematic exploration guided by discrete structure and local optimization. The formulation incorporates torsion angle intervals and chirality constraints through an atom ordering that preserves protein-backbone geometry. Numerical experiments show that the approach efficiently reconstructs geometrically valid conformations even under wide distance bounds, whereas most existing studies assume narrow ones.
	\end{abstract}
	
	\noindent \textbf{keywords:} Distance Geometry, Molecular Distance Geometry Problem, Spectral Projected Gradient Algorithm, 3D Protein Structure, Nuclear Magnetic Resonance
	
	\section{Introduction}\label{sec:dmdgp}
	
	The \emph{Molecular Distance Geometry Problem} (MDGP) aims at the determination of the three-dimensional configuration of a protein molecule consistent with interatomic distance constraints from stereochemistry and Nuclear Magnetic Resonance (NMR) experiments. Typically, only a subset of pairwise distances is available \cite{liberti2014euclidean}.
	
	Formally, let $G=(V,E,d)$ be a weighted graph where $V$ denotes the set of atoms, $E$ the set of (unordered) atom pairs with known distances, and $d:E\to\mathbb{R}_{>0}$ assigns those distances (for an edge $\{u,v\}\in E$ we let $d(\{u,v\})=d_{uv}$). The MDGP asks for a realization $x:V\to\mathbb{R}^3$ such that
	\begin{equation*}
		\|x_u-x_v\|=d_{uv}, \qquad \forall \{u,v\}\in E,
	\end{equation*}
	
	\noindent with $\|\cdot\|$ the Euclidean norm.
	
	A well-studied subclass, the \emph{Discretizable} MDGP (DMDGP), is NP-hard \cite{Lavor2012, liberti2014euclidean} yet exhibits rich combinatorial structure \cite{lavor2021optimality, Liberti2011}. Under an appropriate atom ordering and stereochemical constraints, the DMDGP search space can be represented by a binary search tree that can be explored by the \emph{Branch-and-Prune} (BP) algorithm \cite{Liberti2008}. In this setting, the protein backbone geometry, parameterized by torsion angles $\phi,\ \psi,\ \omega$ \cite{donald2011algorithmsStructMolBios}, is modeled using covalent bond lengths $d_{i-1,i}$ and distances across two bonds $d_{i-2,i}$ given \emph{a priori} \cite{LPB_2021}. From NMR chemical shifts, approximate $\phi,\ \psi$ values can be inferred \cite{Shen2015a} (and typical $\omega$ values are available in structural data \cite{LPB_2021}), which can be converted into distances across three bonds ($d_{i-3,i}$) \cite{Liberti2008}. This completes the model.
	
	To incorporate proximity information from NMR Nuclear Overhauser Effect Spectroscopy (NOESY) (long-range constraints between spatially close hydrogen atoms, $<5$ \AA) \cite{wuthrich1986}, we include backbone-bound hydrogen atoms ($H_N$, $H_\alpha$) together with $N$, $C_\alpha$, and $C$. Because this extension breaks the single-chain structure, we adopt an augmented atom ordering that allows repeated atoms and explicitly includes hydrogen atoms \cite{lavor2019minNMRrigidity,Marques2026}. This ordering also guides BP exploration, since many hydrogen-related constraints span three or more bonds \cite{carvalho2008,Liberti2008}.
	
	Several alternative vertex orders exist \cite{Goncalves2014,liberti2014euclidean}, and in fact we use one of these orders in our experiments. To enhance clarity, however, we will use the standard DMDGP order \cite{Lavor2012} in the theoretical exposition.
	
	Uncertainties in NMR experiments motivate the interval variant (\emph{i}DMDGP) \cite{Lavor2013}, where each $\{i,j\}\in E$ is associated with bounds $[d^L_{i,j},\, d^U_{i,j}]$ and one seeks $x:V\to\mathbb{R}^3$ satisfying
	\begin{equation*}
		d^L_{i,j} \le \|x_i-x_j\| \le d^U_{i,j}, \qquad \forall \{i,j\}\in E.
	\end{equation*}
	
	\noindent In other words, the \emph{i}DMDGP searches for a configuration $X=[x_1,\ldots,x_n]\in\mathbb{R}^{3\times n}$ consistent with interval constraints, which is inherently more challenging than the exact case (from now on, if $d_{ij}^L=d_{ij}^U$ we call $d_{ij}$ an \emph{exact distance}; if $d_{ij}^L<d_{ij}^U$, we call it an \emph{interval distance}).
	
	To illustrate the \textit{i}DMDGP, we consider a toy instance consisting of ten atoms. By definition, the distances $d_{i-1,i}$, $d_{i-2,i}$, $d_{i-3,i}$ are specified for each admissible index $i$ in the DMDGP ordering ($4 \leq i \leq 10$). In this example, we additionally assume that NMR experiments detect interactions between the atom pairs $\{1,8\}$, $\{2,10\}$, and $\{3,9\}$, providing estimates for the corresponding interatomic distances. The distance matrix for this instance is organized as follows: the distances $d_{i-1,i}$ and $d_{i-2,i}$ are exact; all other experimentally determined entries (the “extra distances”) are given as intervals; and unknown entries are denoted by “$?$”, indicating pairs for which no measurement is available. For brevity, we display only the upper triangular part of the symmetric matrix:
	\begin{equation*}
		\left[\begin{array}{cccccccccc}
			0 & d_{1,2} & d_{1,3} & [d_{1,4}^L, d_{1,4}^U] & ? & ? & ? & [d_{1,8}^L, d_{1,8}^U] & ? & ? \\
			& 0 & d_{2,3} & d_{2,4} & [d_{2,5}^L, d_{2,5}^U] & ? & ? & ? & ? & [d_{2,10}^L, d_{2,10}^U] \\
			&  & 0 & d_{3,4} & d_{3,5} & [d_{3,6}^L, d_{3,6}^U] & ? & ? & [d_{3,9}^L, d_{3,9}^U] & ? \\
			&  &  & 0 & d_{4,5} & d_{4,6} & [d_{4,7}^L, d_{4,7}^U] & ? & ? & ? \\
			&  &  &  & 0 & d_{5,6} & d_{5,7} & [d_{5,8}^L, d_{5,8}^U] & ? & ? \\
			&  &  &  &  & 0 & d_{6,7} & d_{6,8} & [d_{6,9}^L, d_{6,9}^U] & ? \\
			&  &  &  &  &  & 0 & d_{7,8} & d_{7,9} & [d_{7,10}^L, d_{7,10}^U] \\
			&  &  &  &  &  &  & 0 & d_{8,9} & d_{8,10} \\
			&  &  &  &  &  &  &  & 0 & d_{9,10} \\
			&  &  &  &  &  &  &  &  & 0 \\
		\end{array}\right].
	\end{equation*}
	
	In most continuous approaches to distance-related problems \cite{Dokmanic2015,Fang2012a,Goncalves2022}, it is assumed that the distance intervals are narrow. In contrast, our framework considers a scenario closer to that observed in NMR experiments, where distance uncertainties are typically wider.
	
	Recent AI-based methods, most notably AlphaFold \cite{alphafold2021}, have profoundly changed protein structure prediction for many folded proteins. However, the problem addressed in this paper is different in nature. We are concerned with the reconstruction of conformations from sparse and uncertain interval restraints, as they arise in NMR-based distance geometry, rather than with sequence-to-single-structure prediction. In this sense, the present work should be viewed as complementary to recent AI-based predictors rather than directly competitive with them.
	
	It is also useful to distinguish between intrinsically disordered proteins and highly flexible proteins \cite{Omidi2026}. The former are more naturally described by ensembles of conformations rather than by a single stable structure, whereas the latter typically possess a well-defined fold but exhibit significant internal motions or transitions between a small number of conformational states. This distinction helps clarify why interval-based distance-geometry formulations remain meaningful in settings involving structural heterogeneity and experimental uncertainty.
	
	In this work, we develop a \emph{hybrid combinatorial--continuous} framework for the \emph{i}DMDGP that integrates the enumerative strategy of \cite{Lavor2013} with a continuous formulation in the spirit of \cite{Glunt1993}. The method couples (i) an enumeration process based on the DMDGP search tree with (ii) a continuous optimization problem whose objective measures violations of admissible intervals. This constitutes a practical step toward scalable hybrid algorithms for realistic \emph{i}DMDGP instances. From a computational viewpoint, we also assess the proposed framework against the interval version of \textsc{MDjeep} \cite{liberti_mdjeep2}, a well-established, publicly available implementation and a representative benchmark for interval distance-geometry methods combining enumeration with continuous refinement.
	
	The rest of this paper is organized as follows. Section~2 presents the problem formulation and the modeling assumptions underlying our hybrid framework. Section~3 details the construction and optimization algorithms, including the angle-guided initialization, sign-consistent improvement, and SPG-based refinement. Section~4 reports the computational experiments, comparing the proposed method with the \textsc{MDjeep} baseline and analyzing runtime behavior and robustness. Finally, Section~5 summarizes the main findings and outlines perspectives for future work.
	
	\section{Problem formulation and modeling assumptions}\label{sec:spg}
	
	We summarize the modeling assumptions and notation used throughout the paper, postponing solver-dependent choices and placement formulas to Sections~\ref{sec:initialXd} and~\ref{sec:tests}.
	
	\subsection{Backbone geometry and atom ordering}
	
	The protein backbone is a known sequence of residues with covalent bond lengths and bond angles, which we fixed to standard values~\cite{donald2011algorithmsStructMolBios,Engh1991}. Variability arises solely from torsion angles. Long-range interatomic distances (typically from NMR) are represented as intervals, whereas distances across one or two covalent bonds are considered exact.
	
	To encode these relations, we adopt an \emph{i}DMDGP ordering with repetitions of atoms, similar to~\cite{lavor2019minNMRrigidity}, which allows us to explicitly include backbone-bound hydrogen atoms\footnote{For proline, the peptide nitrogen is not bound to a hydrogen; we replace $H_N^i$ by a $\delta$-hydrogen attached to $C_\delta$. For glycine, whose $C_\alpha$ is bonded to two hydrogen atoms, only one is retained for modeling purposes.} ($H_N$, $H_\alpha$) and to couple torsion information with distance constraints through consecutive quadruples. This ordering is convenient for representing chirality and other orientation-dependent geometric properties.
	
	\subsection{Torsion-angle domains and chirality}
	
	Let $\tau_i$ denote the torsion angle associated with the $i$th quadruple of consecutive atoms in the ordering. Its feasible set is either a single interval,
	\begin{equation*}
		\mathcal{T}_i=[\tau_i^L,\tau_i^U],
	\end{equation*}
	
	\noindent or the union of two symmetric intervals,
	\begin{equation*}
		\mathcal{T}_i=[-\tau_i^U,-\tau_i^L]\ \cup\ [\tau_i^L,\tau_i^U],
	\end{equation*}
	
	\noindent depending on whether one or both orientations (signs) are admissible. We write $\tau=(\tau_1,\ldots,\tau_m)\in\mathcal{T} = \mathcal{T}_1\times \cdots\times \mathcal{T}_m$ to indicate $\tau_i\in\mathcal{T}_i$ for all $i$. In Section~\ref{sec:initialXd} we exploit the sign structure of $\mathcal{T}_i$ to guide both construction and improvement steps while avoiding erratic side changes.
	
	\subsection{Distance data and mathematical programming formulation}
	
	Let $E$ be the set of pairs of atoms whose distances are known, each endowed with bounds $[d_{ij}^L,d_{ij}^U]$. A conformation $X=[x_1\ldots x_n]\in\mathbb{R}^{3\times n}$ is an \emph{i}DMDGP solution if
	\begin{equation}\label{eq:feas}
		d_{ij}^L \ \le\ \|x_i-x_j\|\ \le\ d_{ij}^U \qquad \text{for all}\quad\{i,j\}\in E.
	\end{equation}
	
	To explicitly represent this experimental uncertainty within the model, we introduce auxiliary variables $d = (d_{ij})_{\{i,j\}\in E}$ defined within the box
	\begin{equation*}
		\Omega \;=\; \big\{\, d\in \mathbb R^{|E|} \ \big|\  d_{ij}^L \le d_{ij} \le d_{ij}^U,\ \{i,j\}\in E \big\}.
	\end{equation*}
	
	\noindent This formulation separates (i) the misfit of $\|x_i-x_j\|$ to the data from (ii) the enforcement of interval bounds. The projection onto $\Omega$ is straightforward, because the projection of $\|x_i-x_j\|$ onto the interval $[d_{ij}^L,d_{ij}^U]$ is simply
	\begin{equation*}
		d_{ij} = \Pi_{[d_{ij}^L,d_{ij}^U]}(\|x_i-x_j\|) = \min\{\max\{\|x_i-x_j\|,d_{ij}^L\},\,d_{ij}^U\}.
	\end{equation*}
	
	The misfit to distance data is quantified by the (nonconvex) \emph{stress} function
	\begin{equation}\label{stress}
		\sigma(X,d) \;=\; \frac{1}{2}\sum_{\{i,j\}\in E} w_{ij}\,\big(\|x_i - x_j\| - d_{ij}\big)^2,
	\end{equation}
	
	\noindent with positive weights $w_{ij}$ (Section~\ref{sec:tests} details our choices).
	
	To determine a three-dimensional protein structure that is consistent with the distance constraints amounts to finding a conformation $X$, and distance variables $d$, such that the stress function $\sigma(X,d)$ is zero and $d$ satisfies the bounds, \emph{i.e.}, $d\in\Omega$. Therefore, our mathematical programming formulation is
	\begin{align}\label{dmdgpcontinuo}
		\min_{X,d}\quad & \sigma(X,d) \\
		\text{s.t.}\quad & d\in\Omega,\quad X\in\mathbb{R}^{3\times n}. \nonumber
	\end{align}
	
	The function $\sigma$ is smooth at any $(X,d)$ with $\lVert x_i-x_j\rVert>0$, for all $\{i,j\}\in E$~\cite{Glunt1993}, which makes it suitable for smooth optimization methods. Nevertheless, computing a global solution of \eqref{dmdgpcontinuo} remains challenging due to its high-dimensional nonconvex landscape. Consequently, practical continuous approaches often focus on generating high-quality initial conformations that can guide local optimization toward near-global minima.
	
	\paragraph{Practical notes.} A solution of \eqref{eq:feas} is invariant under rigid motions. In Section~\ref{sec:initialXd}, we remove this gauge freedom by fixing the first three atoms during the enumerative phase (unlike in the continuous phase, where $X$ remains free in \eqref{dmdgpcontinuo}) and describe how we (i) build angle-aware initial conformations, (ii) apply sign-consistent discrete improvements, and (iii) refine candidates for \eqref{dmdgpcontinuo} via a projected-gradient scheme. Quality metrics -- Largest Distance Error (LDE) and Mean Distance Error (MDE) -- and stopping criteria are reported in Section~\ref{sec:tests}.
	
	\section{Hybrid construction and optimization framework}
	
	\subsection{Background and placement primitive}\label{sec:initialXd}
	
	Generating a near–globally optimal conformation is challenging because $\sigma$ is highly nonconvex. Two families of strategies are typically used to build good initial conformations from \textit{a priori} information.
	
	\smallskip
	
	\noindent \textbf{(i) SDP-based models.} A common approach relies on semidefinite programming (SDP) relaxations for Euclidean distance matrix completion, with notable formulations in \cite{AlHomidan2005,Alfakih1999,Alipanahi2013,Drusvyatskiy2017}. These models aim at finding low-rank realizations that, although lacking a theoretical guarantee, are hopefully three-dimensional. However, they have limitations: (a) generating diverse initial conformations is not straightforward; (b) only distance information is used, without explicit torsion angle signs or chirality; and (c) solving large-scale SDP formulations for the interval case can be computationally demanding. Although facial reduction algorithmic techniques \cite{Drusvyatskiy2017,Huang2017} and low-rank formulations \cite{Bellavia2021a,Burer2003} mitigate costs in the exact case, efficiency remains limited with interval data.
	
	\smallskip
	
\begin{figure}[!htp]
		\centering
		\def\svgwidth{180px}
		\begingroup%
		\makeatletter%
		\providecommand\color[2][]{%
			\errmessage{(Inkscape) Color is used for the text in Inkscape, but the package 'color.sty' is not loaded}%
			\renewcommand\color[2][]{}%
		}%
		\providecommand\transparent[1]{%
			\errmessage{(Inkscape) Transparency is used (non-zero) for the text in Inkscape, but the package 'transparent.sty' is not loaded}%
			\renewcommand\transparent[1]{}%
		}%
		\providecommand\rotatebox[2]{#2}%
		\newcommand*\fsize{\dimexpr\f@size pt\relax}%
		\newcommand*\lineheight[1]{\fontsize{\fsize}{#1\fsize}\selectfont}%
		\ifx\svgwidth\undefined%
		\setlength{\unitlength}{605.00930786bp}%
		\ifx\svgscale\undefined%
		\relax%
		\else%
		\setlength{\unitlength}{\unitlength * \real{\svgscale}}%
		\fi%
		\else%
		\setlength{\unitlength}{\svgwidth}%
		\fi%
		\global\let\svgwidth\undefined%
		\global\let\svgscale\undefined%
		\makeatother%
		\begin{picture}(1,0.77827631)%
			\lineheight{1}%
			\setlength\tabcolsep{0pt}%
			\put(0,0){\includegraphics[width=\unitlength,page=1]{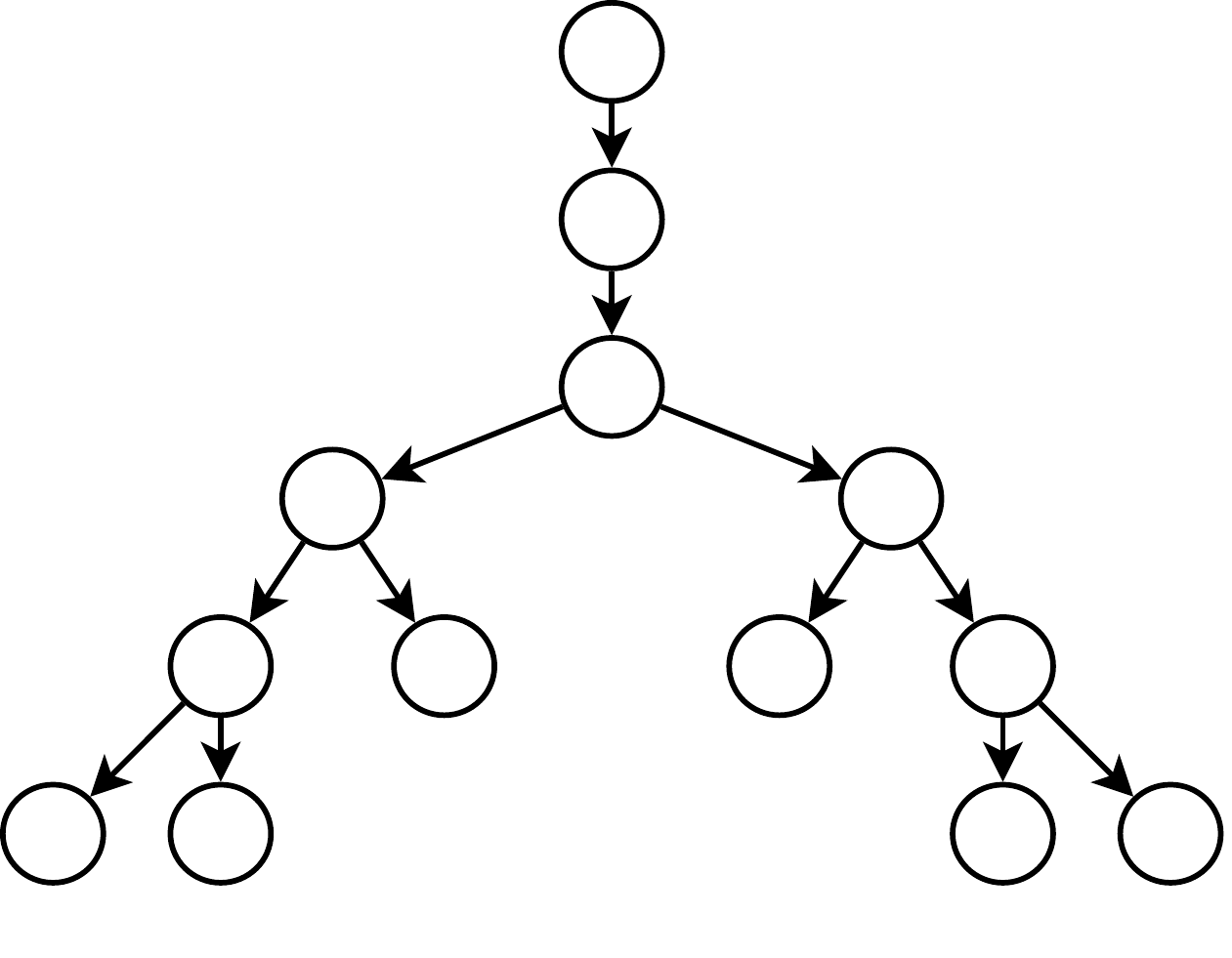}}%
			\put(0.48543685,0.72108078){\color[rgb]{0,0,0}\makebox(0,0)[lt]{\lineheight{1.25}\smash{\begin{tabular}[t]{l}$1$\end{tabular}}}}%
			\put(0.48543685,0.58471924){\color[rgb]{0,0,0}\makebox(0,0)[lt]{\lineheight{1.25}\smash{\begin{tabular}[t]{l}$2$\end{tabular}}}}%
			\put(0.48543685,0.44835768){\color[rgb]{0,0,0}\makebox(0,0)[lt]{\lineheight{1.25}\smash{\begin{tabular}[t]{l}$3$\end{tabular}}}}%
			\put(0.25860816,0.35744434){\color[rgb]{0,0,0}\makebox(0,0)[lt]{\lineheight{1.25}\smash{\begin{tabular}[t]{l}$4$\end{tabular}}}}%
			\put(0.7123202,0.35744434){\color[rgb]{0,0,0}\makebox(0,0)[lt]{\lineheight{1.25}\smash{\begin{tabular}[t]{l}$4$\end{tabular}}}}%
			\put(0.16701492,0.22126223){\color[rgb]{0,0,0}\makebox(0,0)[lt]{\lineheight{1.25}\smash{\begin{tabular}[t]{l}$5$\end{tabular}}}}%
			\put(0.34967335,0.22126223){\color[rgb]{0,0,0}\makebox(0,0)[lt]{\lineheight{1.25}\smash{\begin{tabular}[t]{l}$5$\end{tabular}}}}%
			\put(0.62072696,0.22126223){\color[rgb]{0,0,0}\makebox(0,0)[lt]{\lineheight{1.25}\smash{\begin{tabular}[t]{l}$5$\end{tabular}}}}%
			\put(0.80338542,0.22126223){\color[rgb]{0,0,0}\makebox(0,0)[lt]{\lineheight{1.25}\smash{\begin{tabular}[t]{l}$5$\end{tabular}}}}%
			\put(0.03238519,0.08431139){\color[rgb]{0,0,0}\makebox(0,0)[lt]{\lineheight{1.25}\smash{\begin{tabular}[t]{l}$6$\end{tabular}}}}%
			\put(0.16874672,0.08431139){\color[rgb]{0,0,0}\makebox(0,0)[lt]{\lineheight{1.25}\smash{\begin{tabular}[t]{l}$6$\end{tabular}}}}%
			\put(0.80344776,0.08431139){\color[rgb]{0,0,0}\makebox(0,0)[lt]{\lineheight{1.25}\smash{\begin{tabular}[t]{l}$6$\end{tabular}}}}%
			\put(0.9398093,0.08431139){\color[rgb]{0,0,0}\makebox(0,0)[lt]{\lineheight{1.25}\smash{\begin{tabular}[t]{l}$6$\end{tabular}}}}%
			\put(0.34827652,0.44318879){\color[rgb]{0,0,0}\makebox(0,0)[lt]{\lineheight{1.25}\smash{\begin{tabular}[t]{l}$\tau_4$\end{tabular}}}}%
			\put(0.58381008,0.44318879){\color[rgb]{0,0,0}\makebox(0,0)[lt]{\lineheight{1.25}\smash{\begin{tabular}[t]{l}$-\tau_4$\end{tabular}}}}%
			\put(0.33704033,0.31198258){\color[rgb]{0,0,0}\makebox(0,0)[lt]{\lineheight{1.25}\smash{\begin{tabular}[t]{l}$-\tau_5$\end{tabular}}}}%
			\put(0.14613418,0.31198258){\color[rgb]{0,0,0}\makebox(0,0)[lt]{\lineheight{1.25}\smash{\begin{tabular}[t]{l}$\tau_5$\end{tabular}}}}%
			\put(0.78827304,0.31198258){\color[rgb]{0,0,0}\makebox(0,0)[lt]{\lineheight{1.25}\smash{\begin{tabular}[t]{l}$-\tau_5$\end{tabular}}}}%
			\put(0.59736689,0.31198258){\color[rgb]{0,0,0}\makebox(0,0)[lt]{\lineheight{1.25}\smash{\begin{tabular}[t]{l}$\tau_5$\end{tabular}}}}%
			\put(0.20518879,0.16491378){\color[rgb]{0,0,0}\makebox(0,0)[lt]{\lineheight{1.25}\smash{\begin{tabular}[t]{l}$-\tau_6$\end{tabular}}}}%
			\put(0.03163775,0.17235168){\color[rgb]{0,0,0}\makebox(0,0)[lt]{\lineheight{1.25}\smash{\begin{tabular}[t]{l}$\tau_6$\end{tabular}}}}%
			\put(0.91178955,0.17235168){\color[rgb]{0,0,0}\makebox(0,0)[lt]{\lineheight{1.25}\smash{\begin{tabular}[t]{l}$-\tau_6$\end{tabular}}}}%
			\put(0.7407178,0.16491378){\color[rgb]{0,0,0}\makebox(0,0)[lt]{\lineheight{1.25}\smash{\begin{tabular}[t]{l}$\tau_6$\end{tabular}}}}%
			\put(0,0){\includegraphics[width=\unitlength,page=2]{tree.pdf}}%
		\end{picture}%
		\endgroup%
		\caption{Enumerative search with exact distances. For $i\ge4$, two symmetric placements correspond to the sign of $\tau_i$.}
		\label{fig:BPtree}
	\end{figure}

	\noindent \textbf{(ii) Enumerative methods.} These exploit the discrete nature of the DMDGP by placing atoms sequentially and pruning infeasible branches. The \emph{Branch-and-Prune} (BP) method \cite{Liberti2008} uses only distances; since the quartet $(i-3,i-2,i-1,i)$ does not fix $\tau_i$, two symmetric placements for $x_i$ arise, yielding a binary tree (Figure~\ref{fig:BPtree}). An interval extension (\emph{i}BP) \cite{Lavor2013} samples $d_{i-3,i}$ within $[d_{i-3,i}^L,d_{i-3,i}^U]$ and runs BP for each sample, which can grow combinatorially. To temper branching, \cite{Mucherino2019} integrates local continuous optimization to adjust sampled $d_{i-3,i}$ within their intervals; the \textsc{MDjeep} package \cite{liberti_mdjeep2} applies Spectral Projected Gradient (SPG) method \cite{Birgin2001} to \eqref{dmdgpcontinuo}. Still, angular information is not used explicitly.
	
\begin{figure}[!htp]
		\centering
		\centering
		\begin{subfigure}[b]{0.45\textwidth}
			\centering
			\def\svgwidth{5cm}
			\begingroup%
			\makeatletter%
			\providecommand\color[2][]{%
				\errmessage{(Inkscape) Color is used for the text in Inkscape, but the package 'color.sty' is not loaded}%
				\renewcommand\color[2][]{}%
			}%
			\providecommand\transparent[1]{%
				\errmessage{(Inkscape) Transparency is used (non-zero) for the text in Inkscape, but the package 'transparent.sty' is not loaded}%
				\renewcommand\transparent[1]{}%
			}%
			\providecommand\rotatebox[2]{#2}%
			\newcommand*\fsize{\dimexpr\f@size pt\relax}%
			\newcommand*\lineheight[1]{\fontsize{\fsize}{#1\fsize}\selectfont}%
			\ifx\svgwidth\undefined%
			\setlength{\unitlength}{230.62303162bp}%
			\ifx\svgscale\undefined%
			\relax%
			\else%
			\setlength{\unitlength}{\unitlength * \real{\svgscale}}%
			\fi%
			\else%
			\setlength{\unitlength}{\svgwidth}%
			\fi%
			\global\let\svgwidth\undefined%
			\global\let\svgscale\undefined%
			\makeatother%
			\begin{picture}(1,0.53058622)%
				\lineheight{1}%
				\setlength\tabcolsep{0pt}%
				\put(0,0){\includegraphics[width=\unitlength,page=1]{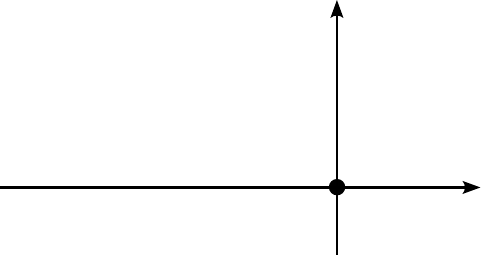}}%
				\put(0.94746151,0.16474707){\color[rgb]{0,0,0}\makebox(0,0)[lt]{\lineheight{1.25}\smash{\begin{tabular}[t]{l}$x$\end{tabular}}}}%
				\put(0.72630528,0.48454794){\color[rgb]{0,0,0}\makebox(0,0)[lt]{\lineheight{1.25}\smash{\begin{tabular}[t]{l}$y$\end{tabular}}}}%
				\put(0.71945737,0.07109221){\color[rgb]{0,0,0}\makebox(0,0)[lt]{\lineheight{1.25}\smash{\begin{tabular}[t]{l}$x_1$\end{tabular}}}}%
				\put(0,0){\includegraphics[width=\unitlength,page=2]{atoms123.pdf}}%
				\put(0.20925536,0.07136109){\color[rgb]{0,0,0}\makebox(0,0)[lt]{\lineheight{1.25}\smash{\begin{tabular}[t]{l}$x_2$\end{tabular}}}}%
				\put(0.42278442,0.16501807){\color[rgb]{0,0,0}\makebox(0,0)[lt]{\lineheight{1.25}\smash{\begin{tabular}[t]{l}$d_{12}$\end{tabular}}}}%
				\put(0,0){\includegraphics[width=\unitlength,page=3]{atoms123.pdf}}%
				\put(0.28544707,0.20243555){\color[rgb]{0,0,0}\makebox(0,0)[lt]{\lineheight{1.25}\smash{\begin{tabular}[t]{l}$\theta_3$\end{tabular}}}}%
				\put(0,0){\includegraphics[width=\unitlength,page=4]{atoms123.pdf}}%
				\put(0.0019524,0.43003491){\color[rgb]{0,0,0}\makebox(0,0)[lt]{\lineheight{1.25}\smash{\begin{tabular}[t]{l}$x_3$\end{tabular}}}}%
				\put(0.01877185,0.22636471){\color[rgb]{0,0,0}\makebox(0,0)[lt]{\lineheight{1.25}\smash{\begin{tabular}[t]{l}$d_{23}$\end{tabular}}}}%
			\end{picture}%
			\endgroup%
			\caption{}
			\label{fig:geometric1}
		\end{subfigure}
		\begin{subfigure}[b]{0.45\textwidth}
			\centering
			\def\svgwidth{7cm}
			\begingroup%
			\makeatletter%
			\providecommand\color[2][]{%
				\errmessage{(Inkscape) Color is used for the text in Inkscape, but the package 'color.sty' is not loaded}%
				\renewcommand\color[2][]{}%
			}%
			\providecommand\transparent[1]{%
				\errmessage{(Inkscape) Transparency is used (non-zero) for the text in Inkscape, but the package 'transparent.sty' is not loaded}%
				\renewcommand\transparent[1]{}%
			}%
			\providecommand\rotatebox[2]{#2}%
			\newcommand*\fsize{\dimexpr\f@size pt\relax}%
			\newcommand*\lineheight[1]{\fontsize{\fsize}{#1\fsize}\selectfont}%
			\ifx\svgwidth\undefined%
			\setlength{\unitlength}{257.05041504bp}%
			\ifx\svgscale\undefined%
			\relax%
			\else%
			\setlength{\unitlength}{\unitlength * \real{\svgscale}}%
			\fi%
			\else%
			\setlength{\unitlength}{\svgwidth}%
			\fi%
			\global\let\svgwidth\undefined%
			\global\let\svgscale\undefined%
			\makeatother%
			\begin{picture}(1,0.54652451)%
				\lineheight{1}%
				\setlength\tabcolsep{0pt}%
				\put(0,0){\includegraphics[width=\unitlength,page=1]{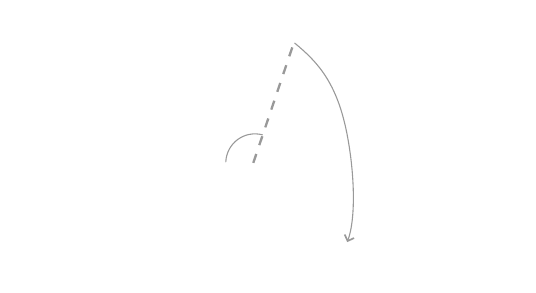}}%
				\put(0.12107915,0.47289316){\color[rgb]{0,0,0}\makebox(0,0)[lt]{\lineheight{1.25}\smash{\begin{tabular}[t]{l}$x_{i-3}$\end{tabular}}}}%
				\put(0,0){\includegraphics[width=\unitlength,page=2]{atomsi.pdf}}%
				\put(0.15032563,0.18585126){\color[rgb]{0,0,0}\makebox(0,0)[lt]{\lineheight{1.25}\smash{\begin{tabular}[t]{l}$x_{i-2}$\end{tabular}}}}%
				\put(0,0){\includegraphics[width=\unitlength,page=3]{atomsi.pdf}}%
				\put(0.50527613,0.2675473){\color[rgb]{0,0,0}\makebox(0,0)[lt]{\lineheight{1.25}\smash{\begin{tabular}[t]{l}$x_{i-1}$\end{tabular}}}}%
				\put(0,0){\includegraphics[width=\unitlength,page=4]{atomsi.pdf}}%
				\put(0.663859,0.04130823){\color[rgb]{0,0,0}\makebox(0,0)[lt]{\lineheight{1.25}\smash{\begin{tabular}[t]{l}$x_i$\end{tabular}}}}%
				\put(0,0){\includegraphics[width=\unitlength,page=5]{atomsi.pdf}}%
				\put(0.4217266,0.14948881){\color[rgb]{0,0,0}\makebox(0,0)[lt]{\lineheight{1.25}\smash{\begin{tabular}[t]{l}$\theta_i$\end{tabular}}}}%
				\put(0.73169636,0.25486053){\color[rgb]{0,0,0}\makebox(0,0)[lt]{\lineheight{1.25}\smash{\begin{tabular}[t]{l}$\tau_i$\end{tabular}}}}%
				\put(0,0){\includegraphics[width=\unitlength,page=6]{atomsi.pdf}}%
				\put(0.94302819,0.4014083){\color[rgb]{0,0,0}\makebox(0,0)[lt]{\lineheight{1.25}\smash{\begin{tabular}[t]{l}$u_1$\end{tabular}}}}%
				\put(0,0){\includegraphics[width=\unitlength,page=7]{atomsi.pdf}}%
				\put(0.8031621,0.35009623){\color[rgb]{0,0,0}\makebox(0,0)[lt]{\lineheight{1.25}\smash{\begin{tabular}[t]{l}$u_2$\end{tabular}}}}%
				\put(0,0){\includegraphics[width=\unitlength,page=8]{atomsi.pdf}}%
				\put(0.76775797,0.48038227){\color[rgb]{0,0,0}\makebox(0,0)[lt]{\lineheight{1.25}\smash{\begin{tabular}[t]{l}$u_3$\end{tabular}}}}%
				\put(0,0){\includegraphics[width=\unitlength,page=9]{atomsi.pdf}}%
				\put(0.21063055,0.28579726){\color[rgb]{0,0,0}\makebox(0,0)[lt]{\lineheight{1.25}\smash{\begin{tabular}[t]{l}$\theta_{i-1}$\end{tabular}}}}%
				\put(0,0){\includegraphics[width=\unitlength,page=10]{atomsi.pdf}}%
			\end{picture}%
			\endgroup%
			\caption{}
			\label{fig:geometric2}
		\end{subfigure}
		\caption{Geometric placement of (a) the first three atoms and (b) the subsequent ones ($x_i$, $i\geq 4$).}
		\label{fig:geometric}
	\end{figure}

	\noindent \textbf{Our stance.} Combinatorial and continuous views are complementary for the \emph{i}DMDGP. We therefore combine (a) fast, angle-aware construction of conformations with (b) continuous refinement by SPG. The geometric placement primitive is as follows. We fix the first three atoms as
	\begin{equation}\label{eqn:x123}
		x_1 = \left[\begin{array}{c}0\\0\\0\end{array}\right], \quad
		x_2 = \left[\begin{array}{c}-d_{12}\\0\\0\end{array}\right], \quad
		x_3 = \left[\begin{array}{c}
			-d_{12} + d_{23}\cos\theta_{3}\\
			d_{23}\sin\theta_{3}\\
			0
		\end{array}\right],
	\end{equation}
	
	\noindent where $\theta_i$ is the bond angle at atom $i$ (between $i$, $i-1$, $i-2$), which satisfies the exact distances among atoms~1–3; see~\cite{Thompson1967}, \cite[proof of Theorem 1]{Lavor2012}. Geometrically, atom~1 is placed at the origin, atom~2 on the negative $x$-axis, and atom~3 in the $xy$-plane, obtained by rotating the segment from $x_2$ to $x_2+(d_{2,3},0,0)$ counterclockwise about $x_2$ by the angle $\theta_3$; see Figure~\ref{fig:geometric}(a). Each subsequent atom $i\ge4$ is placed relative to its predecessors $i-3,i-2,i-1$ as
	\begin{equation}\label{eqn:xi}
		x_i = x_{i-1} + U_i
		\left[\begin{array}{c}
			-d_{i-1,i}\cos\theta_i\\
			d_{i-1,i}\sin\theta_i\cos\tau_i\\
			d_{i-1,i}\sin\theta_i\sin\tau_i
		\end{array}\right],
	\end{equation}
	
	\noindent where the columns of $U_i$ are
	\begin{equation}\label{eqn:xi2}
		u_1 = \frac{v_1}{\|v_1\|}, \qquad
		u_2 = \frac{v_1\times v_2}{\|v_1\times v_2\|}, \qquad
		u_3 = u_2\times u_1,
	\end{equation}
	
	\noindent with $v_1 = x_{i-1} - x_{i-2}$ and $v_2 = x_{i-3} - x_{i-2}$ \cite[expression (7) in Section 4]{Goncalves2014}. The orthogonal matrix $U_i$ represents the change of basis from the canonical basis to basis $\{u_1,u_2,u_3\}$. The vector expanded in~\eqref{eqn:xi} is constructed as follows: we place atoms $i-3$, $i-2$, $i-1$ as before. Atom $i$ is then placed in the same way, based on the positions of atoms~$i-2$ and~$i-1$. Finally, atom~$i$ is rotated about the axis passing through $x_{i-2}$ and $x_{i-1}$ by the torsion angle $\tau_i$.	Multiplication by $U_i$ ``rotates this construction'' to the reference frame defined by the basis $\{u_1,u_2,u_3\}$.	Figure~\ref{fig:geometric}(b) illustrates this construction. After computing $x_i$, we check all long-range distances $d_{hi}$ ($h<i-3$) and prune infeasible branches. This primitive underlies the three algorithms in Sections~\ref{subsec:greedy}–\ref{subsec:multistart}.
	
	\subsection{Algorithms}
	
	Our procedure unfolds in three phases. \emph{Algorithm~\ref{alg:greedy}} builds an angle-aware initial conformation by sampling torsions within their admissible sets and selecting the placement that minimizes a local LDE, defined in~\eqref{localLDE}. \emph{Algorithm~\ref{alg:improve}} performs a discrete neighborhood search by flipping torsion signs whenever both orientations are admissible, restricting subsequent samples to the sign-consistent portion of each torsion domain; only improving flips (in terms of LDE) are accepted. \emph{Algorithm~\ref{alg:multistartspg}} orchestrates multistart runs of Algorithms~\ref{alg:greedy}–\ref{alg:improve}, filters near-duplicates via an RMSD threshold, and refines selected candidates by SPG on the continuous problem \eqref{dmdgpcontinuo} with projected distance variables. Early stopping is triggered as soon as either LDE or MDE (see below) meets the prescribed tolerance.
	
	The choice of SPG over other optimization methods in the refinement phase is motivated by its low computational cost: no system needs to be solved, only first-order derivatives are required, and memory requirements are minimal, making it well suited for large-scale problems such as \eqref{dmdgpcontinuo}.
	Furthermore, it has been shown to effectively minimize general highly nonconvex functions \cite{Birgin2000}, such as \eqref{stress}, up to a moderate accuracy across a wide range of applications (see \cite{Birgin2014,Tavakoli2012} and references therein).
	In particular, SPG was considered in the context of MDGP \cite{Goncalves2022,liberti_mdjeep2,Mucherino2019}.
		
	\subsubsection{Algorithm 1 — Angle-guided greedy construction}\label{subsec:greedy}
	
	Algorithm~\ref{alg:greedy} produces a fast initialization that exploits the chirality information available in the instances of the problem~\cite{Shen2015a}, while remaining consistent with the discretization distances up to numerical precision. For each atom $i$, it samples $\tau_i$ within its admissible domain and selects the trial that yields the smallest local LDE (see expression~\eqref{localLDE} below) with already placed atoms. Since the LDE metric is widely used in the literature as a measure of solution feasibility \cite{2022_bruglieri_MFSoDGC, 2017_dambrosio_NEMaMfRPGfDG, 2018_liberti_BNAiDG}, we recast it here as a local metric that assesses feasibility with respect to the distance information available up to atom $i$.
	\begin{algorithm}[!htp]
		\caption{\textsc{Greedy\_construction}$(N_\text{tors},\mathcal T)$}\label{alg:greedy}
		\smallskip
		\begin{algorithmic}[1]\setstretch{1.1}
			\State Set $x_1$, $x_2$, and $x_3$ as in \eqref{eqn:x123}
			\For {$i=4,\ldots,n$}
			\State $\textsc{Best\_LDE}_i \leftarrow \infty$
			\For {$t=1,\ldots,N_\text{tors}$}
			\State Sample $\tau_i\in \mathcal T_i$ \label{alg:greedysort}
			\State Place $x_i$ using \eqref{eqn:xi}–\eqref{eqn:xi2}
			\If {$\text{LDE}_i(x_1,\ldots,x_i) < \textsc{Best\_LDE}_i$}
			\State $\textsc{Best\_LDE}_i \leftarrow \text{LDE}_i(x_1,\ldots,x_i)$
			\State $\textsc{Best\_torsion}_i \leftarrow \tau_i$
			\EndIf
			\EndFor
			\State $\tau_i\leftarrow \textsc{Best\_torsion}_i$
			\EndFor
			\State Compute $X = [x_1,\ldots,x_n]$ from $\tau$
			\State \Return $\tau$ and $X \in \mathbb{R}^{3\times n}$
		\end{algorithmic}
	\end{algorithm}
	
	That is,
	\begin{equation}\label{localLDE}
		\text{LDE}_i(X)= \max_{\substack{\{i,j\} \in E,\\ j<i}} \delta_{i,j}(X),
	\end{equation}
	
	\noindent where $X$ is a given conformation and, for each edge $\{i,j\}\in E$, the normalized residue $\delta_{i,j}$ is defined as
	\begin{equation*}
		\delta_{i,j}(X) = \max\left\{0,\ \dfrac{d_{ij}^L - \|x_i-x_j\|}{d_{ij}^L},\ \dfrac{\|x_i-x_j\| - d_{ij}^U}{d_{ij}^U}\right\}.
	\end{equation*}
	
	\noindent The quantity above measures the displacement of the computed distance between the atoms $x_i$, $x_j$ and the interval $[d_{ij}^L,d_{ij}^U]$.
	In particular, $\delta_{i,j}(X)=0$ whenever $\|x_i-x_j\|$ lies inside such interval.
	
	\medskip
	
	\noindent\textbf{Global metrics.}
	From $\delta_{i,j}(X)$, we consider the global metrics
	\begin{equation*}
		\mathrm{LDE}(X) := \max_{\{i,j\}\in E}\,\delta_{i,j}(X) \qquad \text{and} \qquad \mathrm{MDE}(X) := \dfrac{1}{|E|}\sum_{\{i,j\}\in E}\delta_{i,j}(X).
	\end{equation*}
	
	For exact distances $\big(d_{ij}^L = d_{ij}^U = d_{ij}\big)$ one has $\delta_{i,j}(X) = \big|\|x_i-x_j\| - d_{ij}\big|/d_{ij}$. The difference between LDE and LDE$_i$ is that the former measures the error along \textit{all} atoms, whereas the latter only takes into account atoms with indices up to $i$, suitable to measure the quality of the partial solutions constructed in Algorithm~\ref{alg:greedy}. In this sense, LDE$_i$ is a  ``local/partial LDE''.
	These quantities are used in the tolerance thresholds $\varepsilon_{\mathrm{LDE}}$ and $\varepsilon_{\mathrm{MDE}}$ in Algorithm~\ref{alg:multistartspg}.
	
	\subsubsection{Algorithm 2 — Sign-consistent local improvement}\label{subsec:improve}
	
	The second phase explores discrete alternatives by flipping torsion signs when both orientations are admissible. To avoid erratic side changes, new trials restrict $\tau_i$ to the portion of $\mathcal T_i$ consistent with the chosen sign:
	\begin{equation*}
		\mathcal T_i^\text{sgn}(\tau) = \begin{cases}
			\mathcal T_i \cap [0,\infty) & \text{if }\tau > 0,\\[0.2cm]
			\mathcal T_i \cap (-\infty,0] & \text{if }\tau < 0.
		\end{cases}
	\end{equation*}
	
	\noindent For each $i$, we flip the sign (if allowed), rebuild a trial conformation with Algorithm~\ref{alg:greedy} using $\mathcal T_i^\text{sgn}$, and keep the change only if it improves the LDE. This targeted move often repairs mismatches introduced by purely local choices in Algorithm~\ref{alg:greedy} without incurring a combinatorial blow-up.
	
	\begin{algorithm}[!htp]
		\caption{\textsc{Improve}$(X, \, \tau, \, N_\text{tors})$}\label{alg:improve}
		\smallskip
		\begin{algorithmic}[1]\setstretch{1.2}
			\For {$i=4,\ldots,n$}
			\State $\tau^\text{trial}\leftarrow \tau$
			\If {$-\tau^\text{trial}_i\in \mathcal T_i$} \Comment{both orientations admissible}
			\State $\tau^\text{trial}_i \leftarrow -\tau^\text{trial}_i$
			\State $\mathcal T_i^\text{trial} \leftarrow \mathcal T_i^\text{sgn}(\tau^\text{trial}_i)$
			\State $\tau^\text{trial}, X^\text{trial} = \textsc{Greedy\_construction}(N_\text{tors},\, \mathcal T_4\times \cdots\times \mathcal T_i^\text{trial}\times \cdots\times \mathcal T_m)$
			\If {$\text{LDE}(X^\text{trial}) < \text{LDE}(X)$}
			\State $X\leftarrow X^\text{trial}$; \ $\tau\leftarrow \tau^\text{trial}$
			\EndIf
			\EndIf
			\EndFor
		\end{algorithmic}
	\end{algorithm}
	
	\subsubsection{Algorithm 3 — Multistart orchestration (Alg.~1+2) and SPG refinement}\label{subsec:multistart}
	
	Algorithm~\ref{alg:multistartspg} diversifies initializations (via multistart of Algorithms~\ref{alg:greedy}–\ref{alg:improve}), discards near-duplicates based on RMSD, and then performs continuous refinement by SPG on \eqref{dmdgpcontinuo}. Distance auxiliaries $d_{ij}$ are initialized by projecting the current interatomic distances onto their admissible intervals $[d^L_{ij},d^U_{ij}]$, yielding a good starting point for the SPG method applied to problem \eqref{dmdgpcontinuo}.
	
	We emphasize that torsion angle information is used to construct informed initial conformations, but it is not imposed as an explicit constraint during the SPG refinement. After Algorithms~\ref{alg:greedy}--\ref{alg:improve} generate a candidate using torsion signs and intervals, the angular intervals are converted into distance intervals that complete the set of distance bounds used in the Cartesian formulation \eqref{dmdgpcontinuo}. This conversion retains the metric information associated with the angular ranges, but not the orientation sign as an explicit constraint. Consequently, the continuous phase may move a conformation outside the original torsion angle ranges. This is a heuristic design choice: angular information reduces the effective branching of the discretized search and provides starting points expected to be close to low-stress conformations, whereas SPG is used to reduce residual distance violations. Thus, torsion signs are not guaranteed to be preserved during refinement.
		
	After each refinement, the candidate is evaluated using the MDE and LDE metrics. The procedure terminates as soon as either tolerance is satisfied; otherwise, the candidate is considered for inclusion in the pool of distinct conformations until the prescribed budget is exhausted.
		
	\paragraph{Similarity filter (RMSD).} We measure similarity between two conformations $X$ and $Y$ via the Kabsch-aligned RMSD:
	\begin{equation*}
		\text{RMSD}(X, Y) = \min_{\substack{Q \in \mathbb{R}^{3\times 3},\\ Q^\top Q = I\phantom{,}}} \dfrac{1}{\sqrt{n}} \big\| \tilde Y - \tilde X Q \big\|_F,
	\end{equation*}
	
	\noindent where $\|\cdot\|_F$ is the Frobenius norm, $Q$ is orthogonal (Kabsch alignment) \cite{Kabsch1976}, and $\tilde X,\ \tilde Y$ are centralized at the same center of mass. For $n\le 200$ we use all atoms; otherwise, only the $C_\alpha$ subset \cite{Aier2016}.
	Candidates with $\text{RMSD}\le \varepsilon_\text{similar}$ are treated as near-duplicates and skipped.
	
	\begin{algorithm}[!htp]
		\caption{Multistart orchestration (Alg.~1+2) with SPG refinement}\label{alg:multistartspg}
		\smallskip
		\noindent\textbf{Initialization.} Set parameters:
		\begin{itemize}
			\item $N_\text{trial} > 0$: maximum number of trials
			\item $N_\text{conf} > 0$: maximum number of distinct initial conformations
			\item $N_\text{tors} > 0$: trials per atom in \textsc{Greedy\_construction}
			\item $N_\text{impr} \geq 0$: repetitions of \textsc{Improve}
			\item $\varepsilon_\text{MDE} > 0$, $\varepsilon_\text{LDE} > 0$:  feasibility  tolerances
			\item $\varepsilon_\text{similar} > 0$: RMSD threshold to declare conformations ``equal''
		\end{itemize}
		
		\noindent Initialize the pool $C\leftarrow\emptyset$
		
		\begin{algorithmic}[1]\setstretch{1.2}
			\For {$c=1,\ldots,N_\text{trial}$}
			\State $\tau^c, X^c = \textsc{Greedy\_construction}(N_\text{tors}, \, \mathcal T)$ \label{alg:multistartspginit} \Comment{Algorithm \ref{alg:greedy}}
			\For {$t=1,\ldots,N_\text{impr}$}
			\State $\textsc{Improve}(X^c, \, \tau^c, \, N_\text{tors})$  \Comment{Algorithm \ref{alg:improve}}
			\EndFor
			\If {$\text{MDE}(X^c) \leq \varepsilon_\text{MDE}$ \textbf{ or } $\text{LDE}(X^c) \leq \varepsilon_\text{LDE}$} \Return $X^c$
			\EndIf
			\If {$\exists\, X^j\in C$ s.t.\ $\text{RMSD}(X^c,\, X^j) \leq \varepsilon_\text{similar}$} \textbf{continue} \label{alg:multistartspgrmsdtest}
			\EndIf
			\State Set $d_{ij}\leftarrow \Pi_{[d_{ij}^L,\, d_{ij}^U]}(\| x_i - x_j \|)$ for each $\{i, j\} \in E$, where $[x_1,\ldots,x_n] = X^c$ \label{alg:multistartspgdij}
			\State Apply SPG to \eqref{dmdgpcontinuo} starting from $(X^c,d)$ \label{alg:multistartspgSPG}
			\If {$\text{MDE}(X^c) \leq \varepsilon_\text{MDE}$ \textbf{ or } $\text{LDE}(X^c) \leq \varepsilon_\text{LDE}$} \Return $X^c$
			\EndIf
			\State $C\leftarrow C\cup \{X^c\}$
			\If {$|C| \geq N_\text{conf}$} \textbf{break}
			\EndIf
			\EndFor
			\State \Return the conformation with the smallest MDE found so far
		\end{algorithmic}
	\end{algorithm}
	
	Since the multistart strategy explores different initialization paths, several distinct conformations can be generated during the search. Each conformation corresponds to a locally optimized structure obtained after refinement. The best conformation found is defined as the one achieving the smallest MDE value among all conformations generated within the computational budget.
	
	\section{Computational results}\label{sec:tests}
	
	\subsection{Experimental setup and parameters}
	
	We implemented all algorithms in \texttt{Julia} and ran the experiments in single-threaded mode on an AMD Threadripper~1950X with 64~GB RAM under Ubuntu~22.04.5 LTS. Source code is available at \url{https://github.com/leonardosecchin/MDGP.jl}.
	
	All weights $w_{ij}$ in \eqref{stress} are set equal, except those associated with discretization distances, which are doubled to preserve, as much as possible, the satisfaction achieved in the construction phase after applying SPG; the weight vector is then normalized.
	
	In Algorithm~\ref{alg:multistartspg}, we declare a ``success'' if, in line~8 of the SPG phase, the stress function value reaches at most $10^{-7}$, or if the stopping criteria in Algorithm~\ref{alg:multistartspg} are met.
	
	We cap SPG at $30{,}000$ iterations per conformation and also adopt a lack-of-progress stop: termination occurs when the step length becomes numerically zero or the relative decrease in stress is insufficient for 100 consecutive iterations. In practice, SPG typically halts well before the iteration cap; hence, ``failure'' corresponds to stalling at a nonoptimal local minimum. SPG parameters follow \cite{Birgin2001}: $\gamma=10^{-4}$, $m=10$, $\lambda_{\min}=10^{-30}$, $\lambda_{\max}=10^{30}$, $\sigma_1=0.1$, $\sigma_2=0.9$; the initial spectral step is computed as in \cite{Birgin2003}.
	
	For Algorithm~\ref{alg:multistartspg} we use
	$N_\text{trial} = 500$,
	$N_\text{conf} = 50$,
	$N_\text{tors} = 20$,
	$N_\text{impr} = 3$,
	$\varepsilon_\text{MDE} = 10^{-3}$,
	$\varepsilon_\text{LDE} = 10^{-2}$, and
	$\varepsilon_\text{similar} = 5.0$.
	We also stop generating new conformations if, for 50 consecutive trials, no conformation different from the previous ones (in the sense of line~\ref{alg:multistartspgrmsdtest}) is found.
	
	For \textsc{MDjeep}, we use the remaining parameters suggested in the example provided in its repository \url{https://github.com/mucherino/mdjeep}, including tolerance $10^{-3}$. To set the resolution parameter, we performed a sensitivity analysis over the values (in~\AA) $0.5$, $1.0$, $2.0$, $3.0$, $4.0$, $5.0$, $6.0$, $8.0$, $10.0$, and $15.0$. The value used in the reported comparison was selected by first maximizing the number of solved instances and, in case of ties, minimizing the total runtime over the solved instances. By this criterion, the selected resolution is $3.0$~\AA.
	
	\subsection{Computational experiments}
	
	We compare Algorithm~\ref{alg:multistartspg} against \textsc{MDjeep} \cite{liberti_mdjeep2}, since both combine enumerative ideas with continuous minimization and \textsc{MDjeep} is a well-established baseline for the setting considered.
	
	The original version of \textsc{MDjeep} \cite{liberti_mdjeep} provides an implementation of the (exact) BP algorithm for the Discretizable DGP (DDGP), a generalization of the DMDGP \cite{liberti_ddgp}. The aim of its continuous component is simply to solve linear systems. In the interval \textsc{MDjeep} \cite{liberti_mdjeep2} the exactness guarantees are weakened by the presence of the intervals, and the continuous component (SPG) is more essential. Therefore, despite the differences of the two algorithms, the analogy between discrete and continuous components is convincing enough to warrant the computational comparison. Since the interval \textsc{MDjeep} makes a stronger effort toward exhaustiveness than ours, one would expect our heuristic to achieve better results than the interval \textsc{MDjeep}, which is exactly what our computational experiments show.
	
	The benchmark set comprises 30 instances, ranging from 87 atoms (17 residues) to 1,347 atoms (269 residues), derived from PDB structures \cite{berman2000ProteinDataBank}. Each instance is constructed by considering angular intervals of width $50^{\circ}$ \cite{malliavin2019} for $\phi$ and $\psi$ torsion angles, and distance constraints of width 1~\AA{} between hydrogen atoms belonging to adjacent residues or 2~\AA{} otherwise \cite{2017_PSbN}. The hydrogen–hydrogen distance constraints are computed exclusively between atoms located within 5~\AA{} in the reference structure.
		
	For reproducibility, each benchmark instance is provided through two plain-text input files, available at \url{https://github.com/leonardosecchin/MDGP.jl}. One file contains the distance data, where each row includes a pair of atom indices $i,j$, the lower and upper distance bounds $d^L_{ij}$ and $d^U_{ij}$. The solver uses only the atom indices and the corresponding bounds. Exact distances are represented by equal lower and upper bounds, whereas interval restraints have $d^L_{ij} < d^U_{ij}$. Another file contains the atom ordering used in the construction phase. Each row specifies the current atom, the three predecessor atoms that form the clique used by the placement formula~\eqref{eqn:xi}, and the torsion information associated with that clique: an orientation sign, a central torsion angle, and an angular tolerance defining a symmetric interval about this central value. When the sign is zero, the corresponding angular information is not enforced. Angular information is also discarded whenever the angular interval, expanded by $20^\circ$ on each side, crosses either $0^\circ$ or $\pm180^\circ$. Thus, the ordering is chosen not only to make the instance discretizable, but also to identify the cliques defining the torsion angles, allowing angular information to be used directly in Algorithms~\ref{alg:greedy}--\ref{alg:improve}.
		
	For computational purposes, the positions of repeated atoms do not need to be explicitly computed. In the proposed implementation, repeated atoms are omitted, and the predecessors of each atom to be placed are taken as the four atoms occupying the three immediately preceding positions in the adopted order. This transformation is therefore only representational, and it leaves the underlying set of feasible geometric solutions unchanged.
	
	Since the benchmark is constructed from PDB-derived instances with reference conformations, it is appropriate for controlled comparisons in terms of feasibility metrics and runtime, but it is not intended as an ensemble-level assessment for intrinsically disordered proteins.
	
	The maximum CPU time for Algorithm~\ref{alg:multistartspg} and \textsc{MDjeep} was set to 5 hours. It is worth noting that runtimes below 1 second were set to 1 second (in such cases, natural fluctuations in measured CPU times make the comparison meaningless). Results are summarized in Table~\ref{tab:results}. Column $N_c$ reports the number of conformations produced by Algorithm~\ref{alg:multistartspg} until termination. When feasibility is not established, each trial conformation is retained provided that it is not similar to those already stored. As soon as a feasible conformation is found, however, the algorithm terminates and declares success. If no feasible conformation is found, the algorithm proceeds until either the maximum number of conformations allowed ($N_{\mathrm{conf}}$) is reached or the limit of $N_{\mathrm{trials}}$ trials is exhausted. All percentages were rounded to the nearest integer to simplify the presentation of results.
	\begin{table}[!htp]
		\small\centering
		\caption{Computational results. Times are in seconds; best values in bold. Italic values for Algorithm~\ref{alg:multistartspg} indicate that the prescribed tolerances were not met; in this case, we report the conformation with the smallest MDE. ``$-$'' indicates that the algorithm does not converge within 5 hours.}
		\label{tab:results}
		\begin{tabular}{lcc|rrrr|rrr}
			& & & \multicolumn{4}{c|}{Algorithm~\ref{alg:multistartspg}} & \multicolumn{3}{c}{\textsc{MDjeep}}\\
			\cline{4-10}
			ID	&	$n$	&	$|E|$	&	$N_c$	&	LDE	&	MDE	&	time	&	LDE	&	MDE	&	time \\
			\hline
			\texttt{1E0Q} & 87 & 363 & 1 & 1.05e-02 & 9.88e-04 & $\bf<1.00$ & 9.90e-03 & 3.30e-04 & $\bf <1.00$\\
			\texttt{1HO7} & 102 & 485 & 1 & 9.50e-03 & 9.60e-04 & $\bf<1.00$ & 6.47e-03 & 2.56e-04 & $\bf<1.00$\\
			\texttt{1LFC} & 127 & 562 & 1 & 1.04e-02 & 9.98e-04 & $\bf<1.00$ & 1.05e-02 & 3.85e-04 & $      115.00$\\
			\texttt{1MMC} & 152 & 702 & 1 & 2.08e-02 & 1.00e-03 & $\bf<1.00$ & 1.00e-02 & 3.96e-04 & $      243.00$\\
			\texttt{6HKA} & 165 & 759 & 1 & 2.35e-02 & 9.86e-04 & $\bf<1.00$ & 2.04e-02 & 4.86e-04 & $        5.00$\\
			\texttt{1SPF} & 177 & 836 & 1 & 3.91e-02 & 9.88e-04 & $\bf<1.00$ & 2.99e-02 & 4.62e-04 & $        4.00$\\
			\hline\texttt{1C56} & 202 & 942 & 3 & 2.16e-02 & 9.99e-04 & $\bf        1.31$ & 1.57e-02 & 4.67e-04 & $       71.00$\\
			\texttt{2LKS} & 220 & 1,017 & 1 & 4.11e-02 & 9.90e-04 & $\bf<1.00$ & 2.49e-02 & 4.65e-04 & $        1.40$\\
			\texttt{7EAU} & 245 & 1,121 & 3 & 2.79e-02 & 9.99e-04 & $\bf        3.77$ & 1.37e-02 & 4.75e-04 & $      698.00$\\
			\texttt{6QBK} & 247 & 1,151 & 7 & 3.82e-02 & 1.00e-03 & $\bf       10.06$ & \text{$-$} & \text{$-$} & $\text{$-$}$\\
			\texttt{2KNX} & 252 & 1,150 & 2 & 1.91e-02 & 1.00e-03 & $\bf        2.12$ & 2.03e-02 & 4.14e-04 & $       87.69$\\
			\texttt{7PQW} & 252 & 1,156 & 14 & 3.71e-02 & 1.00e-03 & $\bf       18.37$ & 1.17e-02 & 4.65e-04 & $     1,830.00$\\
			\hline\texttt{2EGE} & 375 & 1,688 & \it50 & \it2.56e-02 & \it1.00e-03 & $\it      103.80$ & 1.51e-02 & 4.22e-04 & $     6,200.29$\\
			\texttt{2YRT} & 375 & 1,732 & 5 & 3.02e-02 & 1.00e-03 & $\bf        5.13$ & 2.53e-02 & 4.49e-04 & $     7,378.00$\\
			\texttt{2KW9} & 377 & 1,656 & 1 & 5.61e-02 & 9.99e-04 & $\bf<1.00$ & 2.19e-02 & 4.77e-04 & $       17.47$\\
			\texttt{1KCY} & 377 & 1,838 & 6 & 3.09e-02 & 9.99e-04 & $\bf       12.12$ & 6.12e-02 & 4.32e-04 & $      867.59$\\
			\texttt{2N2N} & 475 & 2,302 & 17 & 7.46e-02 & 1.00e-03 & $\bf       37.23$ & \text{$-$} & \text{$-$} & $\text{$-$}$\\
			\texttt{2N9D} & 475 & 2,318 & 5 & 3.31e-02 & 1.00e-03 & $\bf        6.33$ & \text{$-$} & \text{$-$} & $\text{$-$}$\\
			\hline\texttt{2EBT} & 500 & 2,216 & 2 & 4.61e-02 & 1.00e-03 & $\bf<1.00$ & 3.40e-02 & 4.97e-04 & $       29.18$\\
			\texttt{1J0F} & 500 & 2,358 & 30 & 3.28e-02 & 1.00e-03 & $\bf       60.49$ & \text{$-$} & \text{$-$} & $\text{$-$}$\\
			\texttt{8VRC} & 500 & 2,362 & 8 & 2.55e-02 & 9.94e-04 & $\bf       18.98$ & 2.67e-02 & 4.16e-04 & $     2,601.07$\\
			\texttt{2VB5} & 502 & 2,268 & \it50 & \it5.10e-02 & \it1.06e-03 & $\it      149.28$ & \text{$-$} & \text{$-$} & $\text{$-$}$\\
			\texttt{2J4M} & 502 & 2,309 & 5 & 6.45e-02 & 1.00e-03 & $\bf       10.06$ & 1.91e-02 & 4.27e-04 & $     3,354.19$\\
			\texttt{1DX0} & 520 & 2,424 & 1 & 4.08e-02 & 1.00e-03 & $\bf        1.88$ & 3.90e-02 & 5.01e-04 & $     3,902.37$\\
			\hline\texttt{1A66} & 892 & 4,023 & \it50 & \it4.59e-02 & \it1.70e-03 & $\it      214.37$ & \text{$-$} & \text{$-$} & $\text{$-$}$\\
			\texttt{1BC9} & 1,002 & 4,674 & 22 & 6.75e-02 & 1.00e-03 & $\bf      101.90$ & \text{$-$} & \text{$-$} & $\text{$-$}$\\
			\texttt{1JCU} & 1,042 & 4,774 & \it50 & \it5.93e-02 & \it1.71e-03 & $\it      220.41$ & \text{$-$} & \text{$-$} & $\text{$-$}$\\
			\texttt{1AP8} & 1,067 & 5,047 & \it50 & \it7.94e-02 & \it1.93e-03 & $\it      182.61$ & \text{$-$} & \text{$-$} & $\text{$-$}$\\
			\texttt{1EZA} & 1,297 & 6,138 & \it50 & \it4.99e-02 & \it1.29e-03 & $\it      332.25$ & \text{$-$} & \text{$-$} & $\text{$-$}$\\
			\texttt{1AH2} & 1,347 & 6,427 & \it50 & \it8.41e-02 & \it3.65e-03 & $\it      277.80$ & \text{$-$} & \text{$-$} & $\text{$-$}$\\
			\hline
		\end{tabular}
	\end{table}
	
	\paragraph{Runtime behavior.} Despite being much less expensive computationally, Algorithm~\ref{alg:multistartspg} achieves high robustness rates. Within the time limit, \textsc{MDjeep} solves 63\% of the instances, whereas Algorithm~\ref{alg:multistartspg} solves 77\% and does so in markedly smaller runtimes. The runtime of \textsc{MDjeep} varies substantially even for instances of similar size. Enumerative search is highly sensitive to long-range constraints appearing late in the ordering: infeasibilities discovered near the bottom of the tree trigger large reconstructions. Moreover, early inaccurate samples can cascade through the tree, forcing many retries. Our hybrid strategy is less sensitive to these effects: once an initial (possibly infeasible) conformation is built, we optimize it as a whole. Overall runtime depends more on the number of distinct conformations generated before convergence, often small ($\leq 5$), than on the raw number of atoms/distances. This suggests that Algorithms~\ref{alg:greedy} and~\ref{alg:improve} jointly provide strong initializations for SPG.
	
	When interpreting these comparisons, it should be noted that \textsc{MDjeep} and Algorithm~\ref{alg:multistartspg} do not rely on the same input information. Whereas \textsc{MDjeep} performs its search without explicit angular guidance, Algorithm~\ref{alg:multistartspg} exploits torsion angle signs and intervals during the construction phase to generate initial conformations. We nevertheless use \textsc{MDjeep} as the baseline because it is a well-established solver for interval distance geometry that combines discrete search with continuous refinement. The results suggest that angular information can substantially improve this class of hybrid methods, indicating that torsion angle information should be explicitly incorporated into the design of distance geometry solvers for protein instances. 		Indeed, when only interval distance information was available, Algorithm~\ref{alg:multistartspg} failed to solve any of the 16 \textsc{MDjeep} (v0.3) instances with interval data. To enable a meaningful comparison, we reordered the atoms to match the discretization order required by our algorithm and supplemented the instances with torsion angle signs extracted from the corresponding PDB structures. Under these conditions, Algorithm~\ref{alg:multistartspg} successfully solved 14 of the 16 instances (88\%). These results further emphasize the importance of angular information, showing that angle signs are a key ingredient for the effectiveness of our approach.
		
	\paragraph{Role of continuous refinement.} Continuous local optimization is essential. \textsc{MDjeep} consistently handles instances up to $252$ atoms (about 92\% of them), which would be unlikely for purely enumerative methods with wide intervals. Similarly, removing SPG from Algorithm~\ref{alg:multistartspg} (omitting lines~\ref{alg:multistartspgdij}–\ref{alg:multistartspgSPG}) produces a dramatic loss in robustness (indeed, no instance in Table~\ref{tab:results} was solved within the prescribed tolerances only using Algorithms~\ref{alg:greedy} and~\ref{alg:improve}). Even when more distinct initial conformations are generated ($N_\text{conf} = 1,000$) and more improvement trials are allowed ($N_\text{tors}=1,000$ and $N_\text{impr}=50$), the algorithm was able to solve only 3 of the 30 instances. In short, the approach is ineffective without the refinement step via optimization of \eqref{dmdgpcontinuo}.
	
	\subsection{Performance profiles and ablation}
	
\begin{figure}[!htp]
		\centering
		\includegraphics[width=0.65\linewidth]{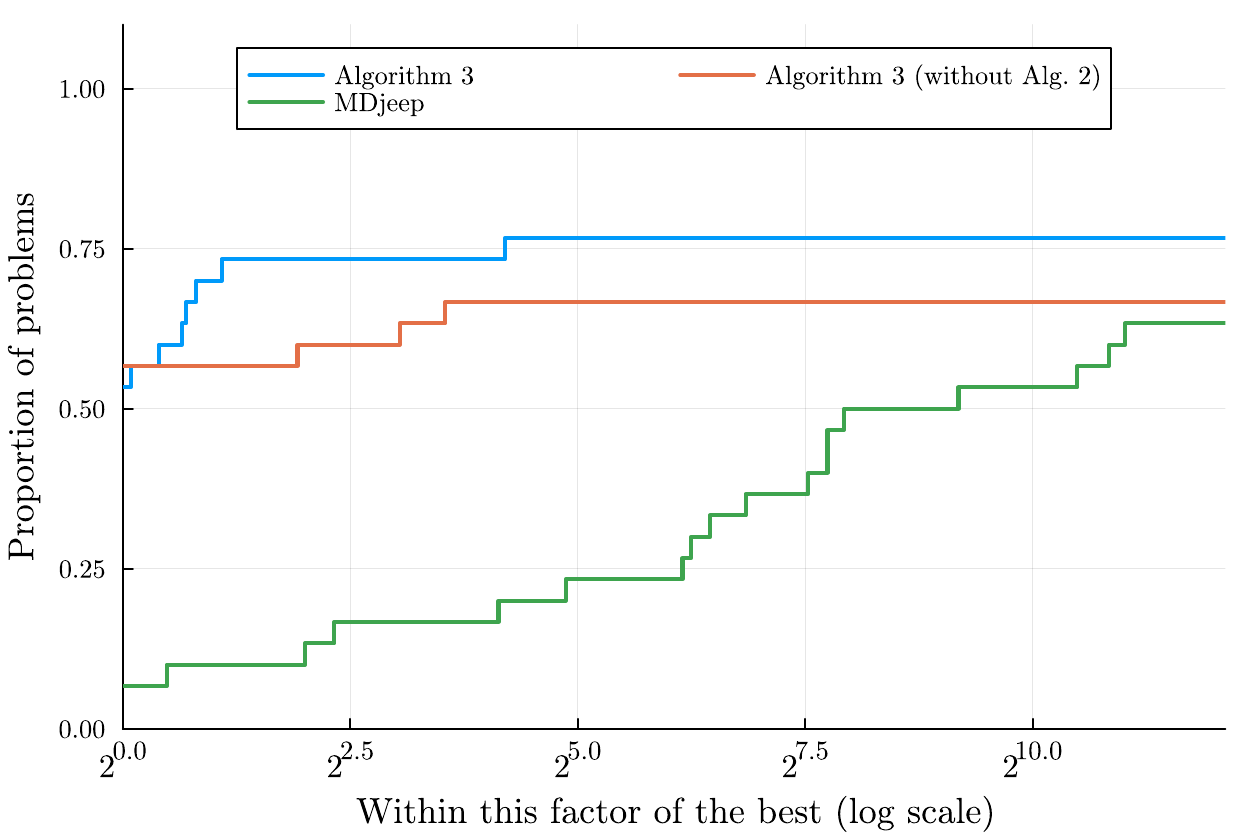}
		\caption{Runtime performance profiles among all instances ($x$-axis in $\log_2$ scale).}
		\label{fig:perfprofile}
	\end{figure}

	We present the results in terms of performance profiles \cite{Dolan2002} (see Figure \ref{fig:perfprofile}), a widely used tool in the optimization community to compare efficiency and robustness of algorithms. We briefly recall their definition. Consider the set $P$ of all instances (in our case, \texttt{1E0Q}, \texttt{1HO7}, \texttt{1LFC}, etc.) and an algorithm $a$. The performance profile of algorithm $a$ is the plot of the function $\rho_a:[1,\infty)\to [0,1]$ defined by
	\begin{equation*}
		\rho_a(t) = \frac{1}{|P|} \, \text{card} \big\{ p\in P\mid r_{p,a}\leq t \big\}, \qquad
		r_{p,a} = \frac{\text{runtime of algorithm $a$ on $p$}}{\text{minimum runtime among all algorithms on $p$}},
	\end{equation*}
	
	\noindent where ``card'' stands for the cardinality of a set. The expression $r_{p,a} = R \geq 1$ means that ``algorithm $a$ takes $R$ times the runtime of the fastest algorithm on $p$''. So, $\rho_a(t)$ is the percentage of problems solved by algorithm $a$ in \emph{at most} $t$ times the runtime of the fastest algorithm.
	
	\paragraph{Two useful reading points.} (i) At $t=1$ (left edge), $\rho_a(1)$ is the fraction of problems on which $a$ is the \emph{fastest}; this quantifies \emph{efficiency}. (ii) For an intermediate value of $t$, for example $t=2^4=16$, $\rho_a(16)$ is the fraction of problems $a$ solves within at most sixteen times the best runtime (note that the $x$-axis in Figure~\ref{fig:perfprofile} is shown on a $\log_2$ scale for better visualization of performance ratios). In this case, \textsc{MDjeep} and Algorithm~\ref{alg:multistartspg} solved $18\%$ and $75\%$ of the problems, respectively. Reading the curve at these points provides complementary views: the left side reflects how often $a$ wins outright; points to the right indicate how tolerant one must be (in runtime multiples) to obtain a given coverage of solved instances. When an algorithm $a$ does not solve $p$, its runtime is defined as $\infty$; then $r_{p,a} = \infty$ (by convention, $r_{p,a} = \infty$ also when ``$\infty/\infty$'' occurs). This is the case of \texttt{1A66} for both algorithms. With this, $\rho_a(t)$ for large $t$ gives the fraction of problems solved by $a$ regardless of runtime (\emph{i.e.}, \emph{robustness}). In this sense, performance profiles balance efficiency (near $t=1$) and robustness (large $t$). It is worth noting that runtimes below 1 second were set to 1 second (as in Table~\ref{tab:results}); in such cases, natural fluctuations in measured CPU times make a finer comparison meaningless. Figure~\ref{fig:perfprofile} shows that, under the 5-hour cap, \textsc{MDjeep} needs more than $2^{11}\!=\!2{,}048$ times the runtime of Algorithm~\ref{alg:multistartspg} to reach its maximum robustness on the instances considered.
	
	\paragraph{Ablation: impact of sign-consistent improvement.} We measure the performance of Algorithm~\ref{alg:multistartspg} without the improvement step (Algorithm~\ref{alg:improve}), which corresponds to setting $N_\text{impr}=0$ (\emph{i.e.}, conformations are computed solely with Algorithm~\ref{alg:greedy}). For a fair comparison, we set $N_\text{tors}=60$ so that the total number of torsion trials is comparable to the regular algorithm. The two variants are similarly efficient (about $53\%$ vs.\ $57\%$ fastest cases, read at $t=1$), but the full algorithm solves more problems overall ($77\%$ vs.\ $67\%$, read on the right side of the profiles). Its performance profile dominates from $t\approx 2^0=1$, indicating that sign-consistent flips improve robustness without harming efficiency.
	
	\section{Final remarks}\label{sec:conclusions}
	
	We proposed a hybrid combinatorial–continuous framework for the \emph{i}DMDGP in protein backbone reconstruction from NMR data. The approach couples (i) an angle-guided greedy construction, (ii) a sign-consistent discrete improvement, and (iii) multistart orchestration with spectral projected gradient (SPG) refinement. By injecting torsion angle interval information into an enumerative placement and then applying continuous optimization, the method balances discrete exploration with local refinement.
	
	\paragraph{Post-AlphaFold context and trade-offs.}
	Recent AI-based predictors, most notably AlphaFold, have transformed protein structure prediction for many folded proteins. Our contribution should be interpreted as complementary to this line of work rather than as a replacement for it. While AlphaFold addresses a sequence-to-single-structure task, the present paper focuses on a restraint-driven inverse problem, namely the reconstruction of geometrically consistent conformations from sparse and uncertain interval distances. This distinction is particularly relevant when experimental uncertainty is substantial and multiple conformations may be compatible with the available data. From a methodological viewpoint, this setting also highlights the trade-offs between continuous and discrete approaches in distance geometry. Early MDGP strategies were primarily based on continuous optimization (see, for example, \cite{Dokmanic2015, Fang2012a, Glunt1993}), whereas the discretizable subclass admits an enumerative treatment \cite{Lavor2012}. Each paradigm has clear strengths and limitations; combining them is therefore natural. Our results support this view: discrete structure sharply narrows the search, while continuous refinement improves agreement with interval data and helps compensate for inaccuracies introduced during the discrete construction phase.
	
	\paragraph{Significance for NMR-based structure determination.} In NMR, distances are indirect and often reflect time-averaged interactions over conformational ensembles rather than a single static structure. NOESY-derived bounds, in particular, are affected by motional averaging and experimental uncertainty, leading to intervals that may be substantially wider than X-ray--based restraints \cite{wuthrich1986}. Methods that remain effective under wide bounds and heterogeneous constraints are therefore useful. Our framework addresses this need by using distance geometry to constrain placements to chemically plausible regions and continuous optimization to reconcile all interval restraints at the conformation level. In practice, the multistart mechanism yields a pool of diverse candidates that can be clustered into conformational sub-states (via RMSD), aligning with ensemble interpretations common in NMR analysis. This also helps explain why such formulations remain relevant alongside powerful AI-based single-structure predictors.
	
	\paragraph{Discrete–continuous synergy.} The experiments highlight a general principle: discrete structure and continuous optimization are complementary for inverse problems in distance geometry. The discretization (ordering, pruning, sign logic) encodes chirality and backbone geometry and reduces the feasible region; the continuous stage provides robustness to noise and exploits all interval information jointly. This division of labor mitigates the combinatorial growth typical of purely enumerative strategies and reduces runtime variability caused by late-violated constraints.
	
	\paragraph{Summary of empirical findings.} On a benchmark of 30 protein instances (87–1,347 atoms), the method solved a larger fraction of problems than the \textsc{MDjeep} baseline under the same 5-hour cap (77\% vs.\ 63\%) and did so with shorter runtimes. Ablation experiments indicate that sign-consistent flips improve robustness without degrading efficiency, and that the SPG stage is necessary to reach embeddings that satisfy LDE/MDE tolerances. These findings support the claim that seeding continuous refinement with torsion-aware discrete construction is effective under wide distance intervals.
	
	\paragraph{Future work.} Our initialization samples torsion angles uniformly, which does not reflect the statistics of real proteins. Because Algorithm~\ref{alg:greedy} is modular (line~\ref{alg:greedysort}), it can incorporate priors from Ramachandran maps or PDB-derived distributions \cite{Lovell2003}. A second direction is to enforce van der Waals lower bounds, either as additional intervals or in a dedicated stage; while this increases the number of stress terms, scalable updates (\emph{e.g.}, stochastic or block-coordinate) \cite{Amaral2022} can limit overhead, or van der Waals checks can be deferred to post-processing. Finally, ensemble-aware formulations that optimize multiple conformations against time-averaged restraints and relaxation-derived constraints would better capture protein dynamics; this can be integrated into our framework by refining a small set of diverse conformers from the multistart pool with shared interval objectives and mild regularization toward diversity. A particularly relevant application domain for such ensemble-aware extensions is that of intrinsically disordered proteins, for which single-structure benchmarks are less suitable and validation should ideally be performed at the ensemble level. Related distance-geometry ideas have also proved useful in protein loop modeling \cite{Labiak2022}.\\
	
	In summary, explicitly combining distance geometry with continuous optimization yields a practical and scalable solver for realistic \emph{i}DMDGP instances, accommodating the wide and uncertain restraints typical of NMR. The hybrid route leverages the complementary advantages of completeness-oriented discrete search and robustness-oriented continuous refinement, and it appears to be a promising direction for distance-geometry problems with experimental uncertainty.
	
	\subsection*{Acknowledgements}
	
	\FUNDING {} We would also like to thank Dr. Marcos Raydan, from the Universidade Nova de Lisboa, for the fruitful discussions about the ideas developed in this paper during CL and MR's visit there in the 1st half of 2025.
	
	\section*{Conflicts of Interest} \par 
	
	The authors declare no conflicts of interest.
		
	\bibliographystyle{abbrvurl}

\end{document}